\documentclass[a4paper,11pt]{article}

\usepackage[german]{babel}
\usepackage{parskip}
\usepackage{blindtext}
\usepackage{amsmath,amsfonts,amssymb,amsthm,epsfig,epstopdf,titling,url,array}
\usepackage{textcomp}
\usepackage{microtype}
\usepackage{graphicx, graphics, color}
\usepackage{wrapfig}
\usepackage{commath}
\usepackage{enumitem}
\usepackage{amsmath}
\usepackage{index}
\usepackage{tikz}
\usepackage{pgfplots}
\usepackage{amssymb}
\usepackage{tkz-euclide}
\usepackage{tikz-dimline}
\usepackage{bm}
\usepackage{subcaption}
\usepackage[dvipsnames]{pstricks} 
\usepackage{pst-solides3d}
\usepackage{pst-3dplot}
\usetikzlibrary{intersections}
\usepackage{cases}
\usepackage{empheq} 
\usepackage{pgfplots}
\usepgfplotslibrary{fillbetween}
\pgfplotsset{compat=1.10}
\usepackage{extarrows}
\usepackage{tikz-3dplot}
\makeindex
\usetikzlibrary{positioning}
\usepackage{algorithm}
\usepackage{caption}
\usepackage{algpseudocode}
\usepackage{gauss}
\usepackage{url}
\usepackage{tabularx}
\usepackage{caption}
\usepackage{lipsum}
\usepackage{titlesec}
\usepackage{cancel}
\usepackage{graphicx} 
\usepackage{caption}  
\usepackage{nicematrix}
\usepackage[section]{placeins}
\usepackage{tabularx}
\usepackage{pdfpages}
\usepackage[a4paper, margin=2.5cm]{geometry}

\usepackage{fancyhdr}
\fancypagestyle{plain}{

\fancyfoot[ER]{\thepage}
\fancyfoot[OR]{\thepage}
}

\newcommand{\R}{\mathbb R}

\newcommand{\headsup}{%
  \mathrel{%
    \vcenter{\offinterlineskip
      \ialign{##\cr$<$\cr\noalign{\kern1pt}$=$\cr}%
    }%
  }%
}
\newcommand{\headdom}{%
  \mathrel{%
    \vcenter{\offinterlineskip
      \ialign{##\cr$>$\cr\noalign{\kern1pt}$=$\cr}%
    }%
  }%
}
\newcommand{\precdom}{%
  \mathrel{%
    \vcenter{\offinterlineskip
      \ialign{##\cr$\prec$\cr\noalign{\kern1pt}$=$\cr}%
    }%
  }%
}
\newcommand{\lessdom}{%
  \mathrel{%
    \vcenter{\offinterlineskip
      \ialign{##\cr$<$\cr\noalign{\kern1pt}$=$\cr}%
    }%
  }%
}

\makeatletter
\newcommand{\thickhline}{
    \noalign {\ifnum 0=`}\fi \hrule height 1pt
    \futurelet \reserved@a \@xhline
}
\newcolumntype{"}{@{\hskip\tabcolsep\vrule width 1pt\hskip\tabcolsep}}
\makeatother
\usepackage[square,sort,comma,numbers]{natbib}
\usetikzlibrary {patterns,patterns.meta}

\title{Berechnung sicherer Fahrradwege}
\author{
  Dr. Julia Sudhoff Santos \\ \small Bergische Universität Wuppertal, Fakultät Mathematik und Naturwissenschaften,\\ \small Gaußstr. 20, 42119 Wuppertal, sudhoff@math.uni-wuppertal.de
  \and
  Lars Kroll \\ \small Bergische Universität Wuppertal, Fakultät Mathematik und Naturwissenschaften, \\
 \small Gaußstr. 20, 42119 Wuppertal, lars.kroll@uni-wuppertal.de
}
\date{}

\begin{document}
\maketitle
\pagestyle{plain}
\section*{Kurzfassung}
Die Sicherheit von Wegen ist nur eingeschränkt messbar und daher schwierig zu quantifizieren. Dahingegen ist es verhältnismäßig leicht Wege bezüglich ihrer Sicherheit in geordnete Kategorien, wie beispielsweise sicher, neutral und gefährlich einzuordnen. 
In diesem Beitrag werden Optimierungsprobleme mit geordneten Kategorien formuliert und Optimalität für diese definiert. Daraus wird eine Lösungsstrategie für solche Probleme abgeleitet. Darüber hinaus wird erklärt, wie die Abgrenzung zwischen den Kategorien erhöht werden kann, sodass längere aber dafür sicherere Wege mit Hilfe von Gewichten berechnet werden können. Diese theoretischen Ergebnisse werden in der Praxis angewendet und es werden auf Grundlage von Daten von OpenStreetMaps sichere Fahrradwege in Stuttgart berechnet. Dabei zeigt sich, dass eine gute Wahl der Gewichte zu weniger Lösungen und kürzeren Rechenzeiten führt.

\section{Motivation und Literaturüberblick}
Es gibt viele Größen, die sich nicht direkt messen lassen, wie beispielsweise die Sicherheit von Wegen für Radfahrer. Es lässt sich zwar erfassen wie viele Unfälle auf bestimmten Straßen passieren, aber dies sagt nur bedingt etwas über die Sicherheit aus. Denn es könnte beispielsweise sein, dass auf einer Straße nur deshalb keine Unfälle passieren, weil sie von Radfahrern nicht genutzt wird. Natürlich kann man die Unfälle in Relation zur Anzahl der Radfahrer setzen, die diese im Schnitt am Tag nutzen und so einen sinnvolleren Messwert ermitteln, aber dann fehlen immer noch Daten über Wege die bisher kaum von Radfahrern genutzt werden. Darüber hinaus hängt die Sicherheit ebenso mit der vorhandenen Infrastruktur zusammen, welche sich schlecht numerisch beschreiben lässt. Dahingegen ist es relativ einfach möglich Wege bezüglich ihrer Sicherheit in geordnete Kategorien einzusortieren. So könnte man entscheiden, dass Hauptstraßen ohne Radweg gefährlich sind, Nebenstraßen ohne Radweg unsicher sind, Straßen mit aufgemalten Radstreifen relativ sicher sind und separate Radwege am sichersten sind. 

Weitere Beispiele für solche geordneten Kategorien, welche auch als ordinale Kosten bezeichnet werden, sind die Klassifizierung der Tierhaltung, der Nutri-Score, olympische Medaillen und Schulnoten. Ordinale Kosten haben die Eigenschaft, dass man nicht weiß um wie viel besser eine Kategorie als eine andere ist. So weiß man nicht wie viel sicherer eine Hauptstraße ohne Radweg verglichen mit einer Nebenstraße ist oder um wie viel besser eine Gold-Medaille verglichen mit einer Silber-Medaille ist, weshalb es auch verschiedene Varianten gibt den Medaillenspiegel bei Olympia aufzustellen. Die meisten Zeitschriften in den Vereinigten Staaten erstellen ein Ranking basierend auf der Gesamtzahl der Medaillen, während das Internationale Olympische Komitee die Länder zunächst nach der Anzahl der Goldmedaillen sortiert. Sollten mehrere Länder die gleiche Anzahl Goldmedaillen haben, werden diese nach den Silbermedaillen sortiert und sollten Länder die gleiche Anzahl Gold- und Silbermedaillen haben, dann werden auch die Bronzemedaillen berücksichtigt. Es gab in der Vergangenheit auch weitere Rankingmethoden, bei denen den Medaillen verschiedene numerische Gewichte zugeordnet wurden (zum Beispiel, dass Goldmedaillen den Wert 5, Silbermedaillen den Wert 3 und Bronzemedaillen der Wert 1 haben), siehe \cite{Olympia1908,Olympia2012}. Wenn den Kategorien numerische Werte zugeordnet werden, welche die Reihenfolge der Kategorien berücksichtigen, spricht man von einer numerischen Repräsentation. Auch Wissenschaftler haben über faire Rankingmethoden nachgedacht, siehe \cite{Sitarz2013medal,GomesJunior2014Sequential,Du2018Modifying,Perini2022Weight}. An diesem Beispiel wird deutlich, dass es für Probleme mit ordinalen Kosten mehr als eine gute Lösung gibt.

Ordinale Kosten wurden bereits in verschiedenen Problemstellungen und mit verschiedenen Definitionen ordinaler Optimalität untersucht. Eine der ersten Untersuchungen ist aus 1971 von Bartee, siehe \cite{Bartee1971Problem}. Danach wurden ordinale Kosten häufig in dem Kontext betrachtet, dass unteilbare Objekte fair auf einige Agenten aufgeteilt werden sollen, siehe \cite{Fishburn1996Binary,Brams2003Fair,Brams2005Efficient,Bouveret2010Fair}. In diesen Referenzen werden verschiedene, meist äquivalente, Definitionen von ordinaler Optimalität verwendet, welche unter anderem auf injektiven Abbildungen und numerischen Repräsentationen beruhen. In anderen Quellen, wie zum Beispiel \cite{Schafer:SP,Klamroth2023Multi}, wird ein Optimalitätskonzept verwendet, welches auf geordneten Vektoren mit den Kategorien in den einzelnen Komponenten basiert. Wird Optimalität auf diesen Vektoren geeignet definiert, so ist dieses Konzept äquivalent zu dem Konzept, welches numerische Repräsentationen nutzt, siehe \cite{OMahony2013Sorted}. Ordinale Kosten wurden auf weitere Probleme angewandt nicht nur auf das bereits beschriebene Zuordnungsproblem. Zum Beispiel werden in \cite{Schweigert1999Ordered} minimale Spannbäume mit ordinalen Kosten berechnet, kürzeste Wege Probleme mit ordinalen Kosten werden in \cite{Bossong1999Minimal,Schafer:SP} untersucht und in \cite{Delort2011Committee,SCHAFER:knapsack} werden unabhängig voneinander ordinale Rucksackprobleme durch eine Umformulierung als multikriterielle Rucksackprobleme gelöst.

Ein Verfahren zum Lösen allgemeiner kombinatorischer Optimierungsprobleme mit ordinalen Kosten durch eine Transformation in ein multikriterielles Optimierungsproblem wird beschrieben in \cite{Klamroth2023Ordinal}. Diese Ergebnisse werden in \cite{Klamroth2023Weighted} verallgemeinert, sodass zusätzliche Informationen über die Gewichtung der Kategorien berücksichtigt werden können. Diese Konzepte aus \cite{Klamroth2023Ordinal,Klamroth2023Weighted} werden im Folgenden beschrieben und auf die Berechnung sicherer Fahrradwege mit beliebig vielen Kategorien angewandt. Das ordinal kürzeste Wegeproblem mit nur zwei Kategorien wurde bereits untersucht und in der Webanwendung \textit{geovelo} umgesetzt, siehe \cite{geovelo,kerg:anef:2021,sauv:sear:2010}.

\section{Optimierung mit ordinalen Kosten}
Es sei ein Optimierungsproblem mit zulässiger Menge $X\subseteq 2^E$ gegeben, welche eine Teilmenge der Potenzmenge einer endlichen disktreten Menge $E$ ist. Jedem Element der Menge $E$ wird eine von $K$ geordneten Kategorien $\mathcal{C}=\{\eta_1,\dots,\eta_K\}$ durch die Abbildung $o:E\to\mathcal{C}$ zugeordnet. Dabei wird angenommen, dass die erste Kategorie die beste ist und jeweils Kategorien mit kleinerem Index besser sind als welche mit größerem Index, das heißt es gilt $\eta_i\prec \eta_{i+1}$ für alle $i=1,\dots,K-1$.

In der Literatur werden mehrere mögliche Zielfunktionen betrachtet, wovon drei im Folgenden skizziert werden. In \cite{Klamroth2023Ordinal} wird gezeigt, dass diese drei Varianten äquivalent sind. Eine Möglichkeit ist es Optimalität über numerische Repräsentationen zu definieren, welche jeder Kategorie ein numerischen Wert zuordnen, sodass die Reihenfolge der Kategorien berücksichtigt wird. Alternativ kann man ein Optimalitätskonzept auf den Vektoren definieren, welche für eine Lösung $x\in X$ zählen, wie viele Elemente dieser Lösung in den verschiedenen Kategorien sind. Im Folgenden betrachten wir ein Konzept, welches die Vektoren $d\in \R^K$ vergleicht, die in der $i$-ten Komponente zählen wie viele Elemente einer Lösung $x=\{e_1,\dots,e_n\}\in X$ in Kategorie $\eta_i$ oder schlechter sind, das heißt $d_i(x)=\vert\{e\in x: o(e)\succeq \eta_i\}\vert$ für alle $i=1,\dots,K$. Auf diesen Vektoren wenden wir Pareto Optimalität an, welche aus der multikriteriellen Optimierung bekannt ist, siehe \cite{Ehrgott2005Multicriteria}. Diese besagt, dass ein Lösungsvektor $d(\hat{x})\in \R^K$ einen anderen Lösungsvektor $d(x')\in \R^K$ dominiert, wenn er in jeder Komponente gleich gut oder besser ist und in mindestens einer Komponente echt besser ist, das heißt $d_i(\hat{x})\leq d_i(x')$ für alle $i=1,\dots,K$ und es existiert mindestens ein Index $j\in\{1,\dots,K\}$ mit $d_j(\hat{x})< d_j(x')$. Ein Lösungsvektor $d(x^*)$ heißt nicht-dominiert, wenn es keinen anderen Lösungsvektor $d(\hat{x})$ gibt, der $d(x^*)$ dominiert. Die zugehörige Lösung $x^*\in X$ heißt effizient. Zusammengefasst ergibt sich folgendes ordinale Optimierungsproblem
\begin{equation}\label{eq:OOP}
	\begin{array}{rl}
		\min & d(x)\\
		\text{s.\,t.} & x\in X,
	\end{array}
\end{equation}
von dem wir alle nicht-dominierten Lösungsvektoren mit einer zugehörigen effizienten Lösung berechnen wollen. Wir erhalten eine Menge an Lösungen, welche untereinander nicht vergleichbar sind.

Im Praxisfall der sicheren Radwege, ist $E$ die Menge von Straßenabschnitten und $X$ die Menge der Routen von Start zum Ziel. Den Straßenabschnitten werden neben den geordneten Kategorien auch noch Längen mittels der Abbildung $l:E\to \R_>$ zugeordnet. Dann zählt der Vektor $d(x)$ nicht die Anzahl der Straßenabschnitte einer Lösung $x=\{e_1,\dots,e_n\}\in X$ die in einer bestimmten Kategorie oder schlechter sind, sondern es werden die Längen dieser Straßenabschnitten zusammen gerechnet, das heißt
\begin{equation}
    d_i(x)=\sum_{e\in x\text{ mit } o(e)\succeq \eta_i} l(e)
\end{equation}
für alle $i=1,\dots,K$. Äquivalent lässt sich der Lösungsvektor auch schreiben als \begin{equation}
d(x)= \sum_{e\in x} l(e)\cdot k(e) \text{ mit } k_i(e)=\begin{cases}
        1, &\text{ wenn } o(e)\succeq \eta_i\\
        0, &\text{ sonst}
    \end{cases} \ \text{ für } i=1,\dots,K.
\end{equation}
Der Vektor $k(e)\in\R^K$ enthält folglich die Information, in welcher Kategorie der Abschnitt $e\in E$ eingeordnet ist. Die resultierenden Lösungsvektoren können dann  mittels der Pareto Optimalität miteinander verglichen werden.

Bisher ist nur bekannt, dass eine Kategorie mit kleinerem Index echt besser ist als eine Kategorie mit größerem Index. Wir werden sehen, dass dies in der Praxis zu einer großen Anzahl an effizienten Lösungen führen kann, die nicht alle eine geeignete Lösung repräsentieren.  So gehört der kürzeste Weg immer zu den sicheren Wegen. Aber in der Praxis möchte man häufig eine klare Abgrenzung zwischen den Kategorien und ist immer bereit kleine Umwege über bessere Straßen in Kauf zu nehmen, um sehr unsichere Straßen zu vermeiden. Dies gilt insbesondere wenn man mit Kinder unterwegs ist.
Es ist daher realistisch Grenzwerte, im folgenden als Gewichte bezeichnet, anzugeben, bis zu denen man immer bereit ist, Umwege über Straßen der nächste besseren Kategorie in Kauf zu nehmen. Dies schränkt die Menge und damit die Anzahl der Lösungen auf eine sinnvolle Teilmenge ein.

Formal nehmen wir an, dass Gewichte
$\omega_i\geq 1$ für $i=1,\dots,K-1$ gegeben sind.  
Das Gewicht $\omega_i\geq 1$, $i\in\{1,\dots,K-1\}$ gibt an um welchen Faktor ein Umweg, der nur über Wegabschnitte der Kategorie $\eta_i$ führt, höchstens länger sein darf, damit er über einen Wegabschnitt der Kategorie $\eta_{i+1}$ bevorzugt wird. 
Angenommen es gibt nur zwei Kategorien, also sichere und unsichere Wege und folgenden drei Routen: Route $\hat{x}$ führt über 6 Kilometer sichere Straße und 1 Kilometer unsichere Straße und Route $\tilde{x}$ führt nur über 8 Kilometer sichere Straße. Gilt $\omega=1$ so sind beide Routen nicht vergleichbar. Hingegen gilt für
$\omega=2{,}3$, dass nun der Weg $\tilde{x}$ besser ist als der Weg $\hat{x}$, weil der Umweg um einen Kilometer unsichere Straße zu vermeiden nur doppelt so lang ist. 

Unter Berücksichtigung dieser zusätzlichen Informationen ergibt sich folgender gewichteter Lösungsvektor $d^{(\omega)}(x)$ für eine Lösung $x=\{e_1,\dots,e_n\}\in X$, siehe \cite{Klamroth2023Weighted}:
\begin{equation}
    d^{(\omega)}_i(x)=\sum_{e\in x} l(e)\cdot k^{(\omega)}(e) \text{ mit } k^{(\omega)}_i(e)=\begin{cases}
        \prod_{\ell=i}^{j-1} \omega_\ell, &\text{ wenn } o(e)= \eta_j,\ i\leq j\\
        1, &\text{ wenn } o(e)= \eta_j,\ i=j\\
        0, &\text{ wenn } o(e)= \eta_j,\ i>j.
    \end{cases}
\end{equation}
Auf diesen Lösungsvektor wird wiederum Pareto Optimalität angewendet. Wenn gilt, dass $\omega_i=1$ für $i=1,\dots,K-1$ erhält man den gleichen Vektor $d$ wie zuvor ohne Gewichte. Die Lösungsvektoren sind für die Beispielrouten in Tabelle~\ref{tab:my_label} angegeben. Für genauere Informationen zu den gewichteten Lösungsvektoren wird auf \cite{Klamroth2023Weighted} verwiesen. 

\begin{table}[H]
    \centering
    \begin{tabular}{c|c|c}
         Route & $\hat{x}$ & $\tilde{x}$\\
         \hline
         $d^{(1)}$ & $\begin{pmatrix}
             7 \\ 1
         \end{pmatrix}$ & $\begin{pmatrix}
             8 \\ 0
         \end{pmatrix}$\\[3ex]
         $d^{(2{,}3)}$ &  $\begin{pmatrix}
             8{,}3 \\ 1
         \end{pmatrix}$& $\begin{pmatrix}
             8 \\ 0
         \end{pmatrix}$
    \end{tabular}
    \caption{Lösungsvektoren für zwei Routen mit zwei Kategorien zu verschiedenen Gewichten~$\omega$.}
    \label{tab:my_label}
\end{table}

Insgesamt erhalten wir folgendes gewichtete ordinale Optimierungsproblem
\begin{equation}\label{eq:wOOP}
	\begin{array}{rl}
		\min & d^{(\omega)}\\
		\text{s.\,t.} & x\in X.
	\end{array}
\end{equation}
Die Optimierungsprobleme~\eqref{eq:OOP} und \eqref{eq:wOOP} lassen sich mit bekannten Verfahren der multikriteriellen Optimierung lösen. Wir lösen die sichersten Wegeprobleme mit dem Multikriteriellen Dijkstra Algorithmus, welcher in \cite{Casas2021Improved} beschrieben wird.

\section{Datengrundlage}\label{sec:Daten}

Um die oben beschriebenen ordinalen kürzesten Wege zu berechnen, benötigen wir Informationen über die Länge und Sicherheit von Straßen. Zu diesem Zweck haben wir die Daten von \emph{OpenStreetMap} \cite{OpenStreetMap} (OSM) verwendet. OSM ist eine Open-Source Datenbank mit Open Data Commons Open Database License (ODbL), die mit Hilfe von vielen Freiwilligen aus vielen Ländern eine große Menge an Geodaten zur Verfügung stellt. 

Mit Hilfe der Python-Bibliothek \emph{OSMnx} \cite{osmnx} ist es möglich, Graphen von Städten oder um einen bestimmten Punkt herum auszugeben. Dieser Graph enthält dann alle Informationen, die auch OpenStreetMaps enthält. Da wir uns für die Berechnung sicherer Fahrradwege interessieren, haben wir den Netzwerktyp ''Fahrrad" gewählt, welcher alle Straßen herausfiltert, die nicht mit dem Fahrrad befahrbar sind. Dies gilt z.B. für Wege mit Treppen sowie für Autobahnen.

Durch die Auswahl des Netzwerktyps „Fahrrad“ werden bereits einige Kanten herausgefiltert, jedoch kann ein Graph einer ganzen Stadt sehr schnell enorme Ausmaße annehmen, z.B. hat der Graph der Stadt Stuttgart (mit Netzwerktyp „Fahrrad“) insgesamt 313.400 Kanten. Für Graphen dieser Größenordnung ist die Berechnung von Routen möglich, jedoch mit längeren Rechenzeiten verbunden. So benötigt beispielsweise die Routenberechnung vom Stuttgarter Hauptbahnhof zum Haus der Wirtschaft in Stuttgart 514 Sekunden. Wenn wir den Graphen hingegen mit einer Funktion der \emph{OSMnx} Bibliothek vereinfachen, erhalten wir einen Graphen mit nur 95.186 Kanten und für die Berechnung der gleichen Route vom Stuttgarter Hauptbahnhof zum Haus der Wirtschaft in Stuttgart werden nur 165 Sekunden benötigt. Der Unterschied in der Rechenzeit ist bei längeren Routen noch signifikanter. Daher nutzen wir vereinfachte Graphen.

Bei der Vereinfachung des Graphen werden Knoten, die keine Kreuzungen oder Sackgassen sind, entfernt und die angrenzenden Kanten dieser Knoten zusammengeführt. Das Zusammenführen erfolgt durch die Erzeugung einer neuen Kante, die zwei benachbarte Kanten eines entfernten Knotens ersetzt. Dabei gehen keine potentiellen Wege verloren.

Durch den Wegfall der Knotenpunkte verlieren wir leider etwas an Genauigkeit bei der Wahl der Start- und Zielknoten, da immer der nächstgelegene Knotenpunkt zur gewählten Start- bzw. Zieladresse gewählt wird. Dadurch kann es vorkommen, dass der Start- bzw. Zielknoten weiter von der tatsächlich gewählten Adresse entfernt liegt als bei nicht vereinfachten Graphen. Die Entfernungen über die wir sprechen liegen in der Regel unter 50 Metern. Ein weiterer Nachteil ist, dass die Kategorisierung unter Umständen nicht so genau ist wie es bei einem nicht vereinfachten Graphen der Fall wäre. Wenn zwei Kanten aus zwei verschiedenen Kategorien zusammengeführt werden, kann der neuen Kante nur eine Kategorie zugeordnet werden. Auf die genaue Art der Kategorisierung wird im nächsten Kapitel näher eingegangen.

Zusätzlich ist noch zu erwähnen, dass es durchaus vorkommen kann, dass es zwischen zwei Knoten mehrere Kanten gibt. Diese werden so gefiltert, dass es zwischen zwei Knoten maximal eine Kante gibt, da die verwendeten Algorithmen nur eine Kante zwischen einem Knotenpaar erwarten. Wir löschen alle Kanten zwischen zwei Knoten, bis auf die Kante in der besten Kategorie (sollte es mehrere geben, wird davon die kürzeste gewählt). Es kann vorkommen, dass zwischen einem Knotenpaar mehrere ordinal nicht-dominierte Kanten liegen. Neben der Kante in der besten Kategorie, könnte es eine kürzere Kante in einer schlechteren Kategorie geben. Diese Option streichen wir durch das Löschen der Kante. Um das Löschen zu verhindern und trotzdem keine Doppelkanten im Graphen zu haben, müsste man den Graphen modifizieren und einen künstlichen Knoten auf den zusätzlichen Kanten einfügen. Da der Aufwand relativ hoch ist, ohne die Ergebnisse signifikant zu verbessern, haben wir darauf verzichtet.

Um die Betrachtung überflüssiger Knoten und Kanten zu vermeiden, ist es meist nicht sinnvoll die Routen auf dem gesamten Graphen einer Stadt zu berechnen, beispielsweise wenn die Route zwischen zwei nah beieinander liegenden Punkten berechnet werden soll. Daher wird in diesem Artikel immer nur der Graph betrachtet, der sich aus einer Entfernung und der Startadresse ergibt. Der entstehende Graph ergibt sich, indem man sich vom Startpunkt aus mit der entsprechenden Distanz in die Richtungen Norden, Süden, Osten und Westen bewegt. Der resultierende Bereich bildet ein Quadrat, das den Graphen repräsentiert, den wir analysieren werden. Das Quadrat ist so zusammengesetzt, dass die Endpunkte der zurückgelegten Strecke in jeder gegebenen Richtung (Norden, Osten, Süden und Westen) den Mittelpunkten der jeweiligen Seiten des Quadrats entsprechen. Dementsprechend werden wir im Folgenden die Größe des betrachteten Graphen mit einer Distanz in Metern und den entsprechenden Startpunkt angeben.

\begin{figure}[!htbp]
\centering
\begin{subfigure}{.5\textwidth}
  \centering
  \includegraphics[width=0.5\textwidth]{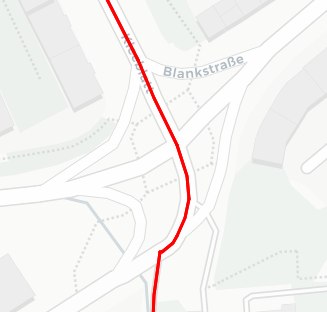}
  
  \caption{Mögliche Route}
  \label{fig:sub1}
\end{subfigure}%
\begin{subfigure}{.5\textwidth}
  \centering
  \includegraphics[width=0.5\textwidth]{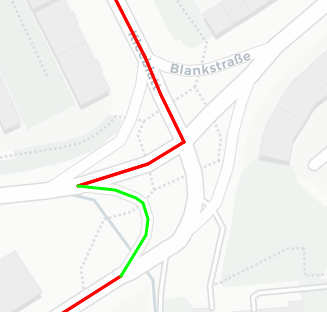}
  
  \caption{Unmögliche Route}
  \label{fig:sub2}
\end{subfigure}

\caption{Zwei verschiedene Routen führen über die gleiche Kreuzung, wobei die rechte in der Praxis nicht nutzbar ist. (Map data © OpenStreetMap contributors \cite{OpenStreetMap})}
\label{fig:Kreuzungen}
\end{figure}

Die Daten von OSM sind sehr hilfreich und machen dieses Projekt erst möglich, allerdings muss auch erwähnt werden, dass die Daten teilweise nicht optimal für unsere Zwecke sind. So sind beispielsweise einige Kreuzungen nicht ideal modelliert, was zum Teil Routen zur Folge hat, die nicht sinnvoll sind. Es kommt zum Beispiel an großen Kreuzungen mit separaten Abbiegespuren vor, dass die Route die Geradeaus-Spur wählt und dann mitten auf der Kreuzung doch rechts abbiegt, siehe Abbildung~\ref{fig:Kreuzungen}. Solche Routen sind nicht nutzbar in der Praxis, aber es ist leicht möglich sie zu einer praxistauglichen Route abzuändern.

\section{Kategorisieren von Kanten}

Nachdem wir nun die Quelle der Graphen und der zugehörigen Daten identifiziert haben, befassen wir uns mit der Kategorisierung der einzelnen Kanten. Ziel der Kategorien ist es die Sicherheit der Kanten aus Sicht eines Radfahrers zu bewerten und in vier Kategorien von sehr gut bis sehr schlecht einzuteilen. Um die Sicherheit einer Straße zu bewerten, verwenden wir die Daten auf den Kanten des OSM Graphen, insbesondere die Information, ob ein Radweg vorhanden ist. In OSM Graphen gibt es unter anderem folgende Attribute „bicycle“, „cycleway”, „cycleway left”, „cycleway right”, „bicycle road”, „cycleway right bicycle” und „cycleway left bicycle”. Wenn eines dieser Attribute einen anderen Wert als „none“ hat, folgt daraus, dass die entsprechende Straße irgendeine Art von Radweg hat. Dadurch können wir die Kanten in zwei Kategorien einteilen, zum einen in Kanten mit Radweg und zum anderen in Kanten ohne Radweg.

Diese beiden Kategorien werden nochmals unterteilt. Bei den Radwegen unterscheiden wir zwischen einem normalen Radweg und einem Radweg, der vollständig vom normalen Verkehr getrennt ist. Ein Radweg ist vom Verkehr getrennt, wenn die Attribute „bicycle“ und „bicycle road” den Wert „yes“ oder „designated“ enthalten.

Die Kanten ohne Radweg werden bezüglich des Straßentyps kategorisiert, welcher im Attribut „highway“ in OSM gespeichert ist. Es gibt z.B. die Attribute „motorway“ (Autobahn), „residential“ (Wohnstraße), „tertiary“ (Landstraße) und viele mehr. Wir haben uns dafür entschieden die Werte „residential“ (Straße, mit der Hauptfunktion Zugang zu Wohngebäuden zu ermöglichen), „living street“ (eine Straße, auf der Fußgänger Vorrang vor Autos haben), „track“ (Waldstraße oder Trassen auf denen Autos kein Zugang haben) und „bridleway“ (ebenfalls eine Art Waldstraße ohne Autos) besser einzustufen als andere Werte. Das bedeutet das Kanten ohne Radweg, aber mit dem Attribut „residential“, „living street“, „bridleway“ oder „track“ in eine bessere Kategorie eingestuft werden als Straßen ohne Radweg und ohne eines dieser Attribute.

Insgesamt ergeben sich dadurch vier Kategorien, wobei die beste Kategorie alle Kanten enthält, die über einen vollständig vom Verkehr getrennten Radweg verfügen, die zweitbeste Kategorie alle Kanten enthält, die über einen Radweg verfügen, aber nicht die Bedingungen der ersten Kategorie erfüllen. Die dritte Kategorie umfasst alle Kanten, die keinen Radweg haben, aber eine „sichere“ Straße sind, und die letzte Kategorie umfasst alle Kanten, die noch keiner Kategorie zugeordnet wurden. 

\begin{figure}[!htbp]
\centering
\begin{subfigure}{.5\textwidth}
  \centering
  \includegraphics[width=0.7\textwidth]{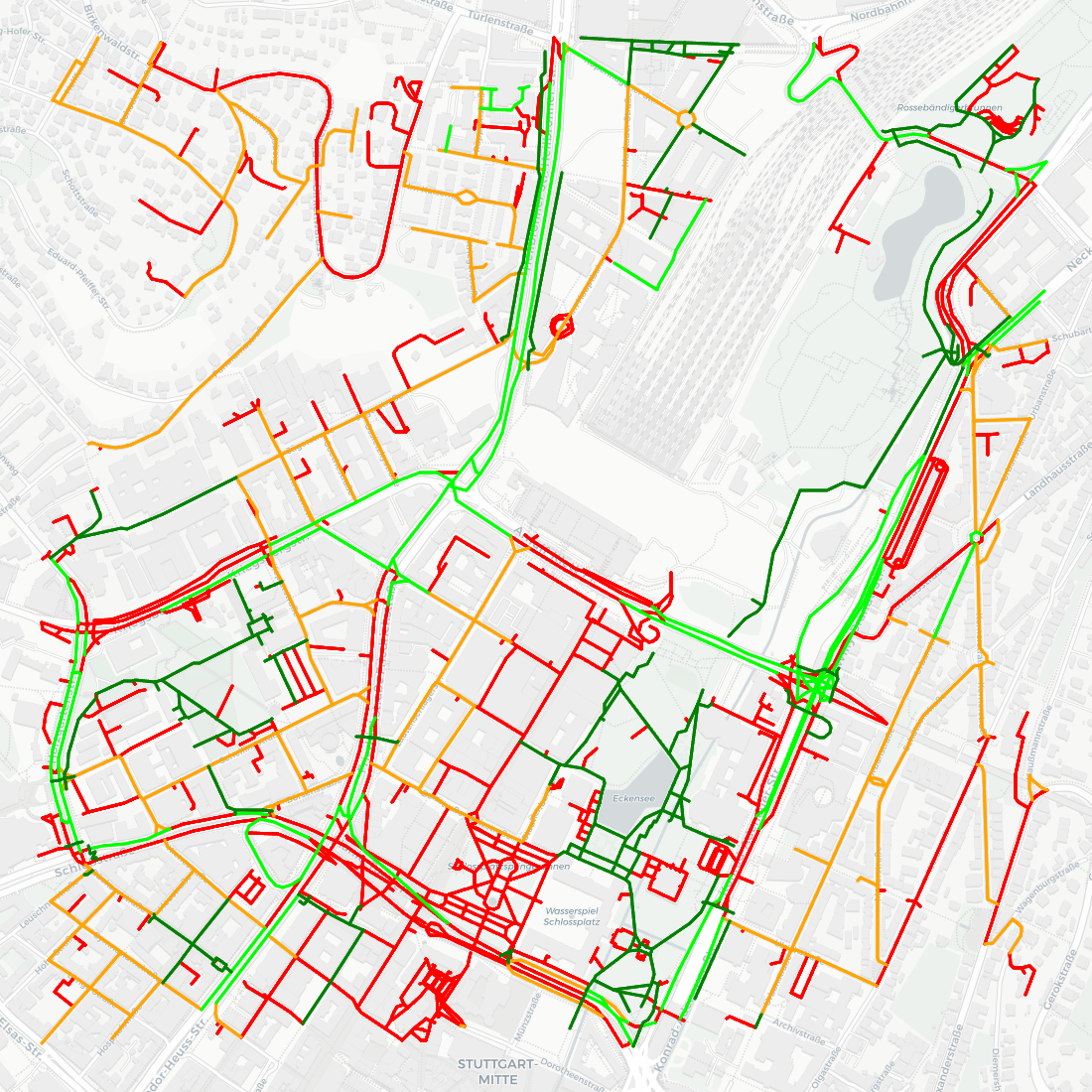}
  
  \caption{800 Meter }
  \label{fig:stuttgart1}
\end{subfigure}%
\begin{subfigure}{.5\textwidth}
  \centering
  \includegraphics[width=0.7\textwidth]{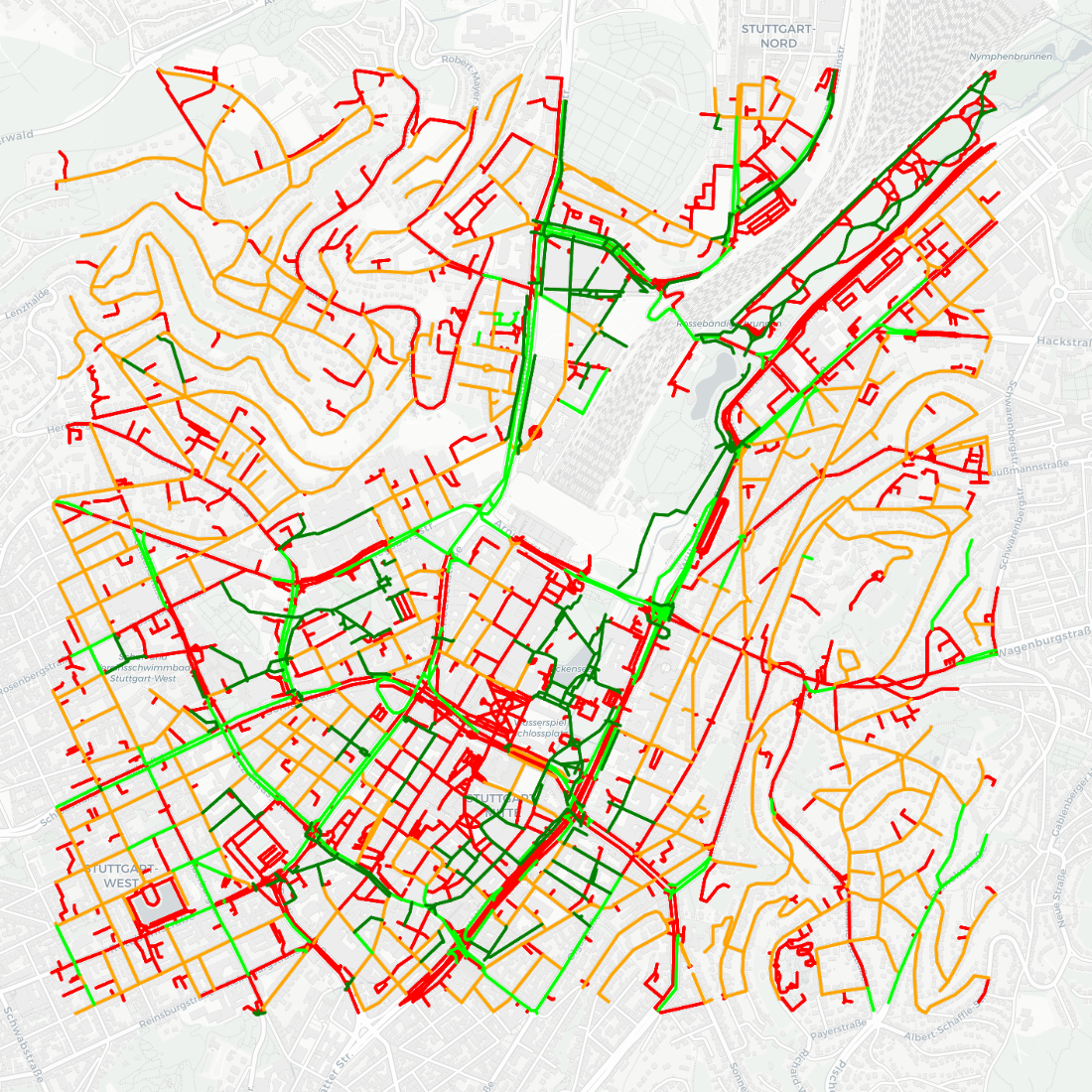}
  
  \caption{1.500 Meter}
  \label{fig:stuttgart2}
\end{subfigure}

\caption{Der Graph von Stuttgart mit dem Hauptbahnhof Eingang im Zentrum und der entsprechenden Graphen Größe. Alle Kanten sind der jeweiligen Kategorie zugeordnet und entsprechend farblich markiert. (Map data © OpenStreetMap contributors \cite{OpenStreetMap}) }
\label{fig:stuttgart}
\end{figure}

Die Abbildung~\ref{fig:stuttgart} zeigt den vollständig kategorisierten Graphen von Stuttgart mit dem Eingang zum Hauptbahnhof als Zentrum in verschiedenen Größen. Alle Kanten sind entsprechend ihrer Kategorie farblich markiert, wobei die Farben die Kategorien in folgender Reihenfolge von sehr gut bis sehr schlecht repräsentieren: dunkelgrün, hellgrün, orange und rot.

\section{Berechnung von Routen}

In diesem Abschnitt werden ungewichtete und gewichtete ordinal effiziente Routen vom Stuttgarter Hauptbahnhof zum Haus der Wirtschaft (dem Veranstaltungsort der Heureka 2024) berechnet. Als Startpunkt wird die „Arnulf-Klett-Passage, 70173 Stuttgart“ (der Haupteingang des Hauptbahnhofs) und als Zielpunkt die „Willi-Bleicher-Straße 19, 70174 Stuttgart“ gewählt. Die Luftlinie zwischen diesen beiden Punkten beträgt ca. 0,78 km. Die Route zwischen diesen beiden Punkten wird für einen Graphen der Größe 800 und 1.500 Meter berechnet, siehe Abschnitt~\ref{sec:Daten}.

\begin{figure}[!htbp]
\centering
\begin{subfigure}{.5\textwidth}
  \centering
  \includegraphics[width=0.78\textwidth]{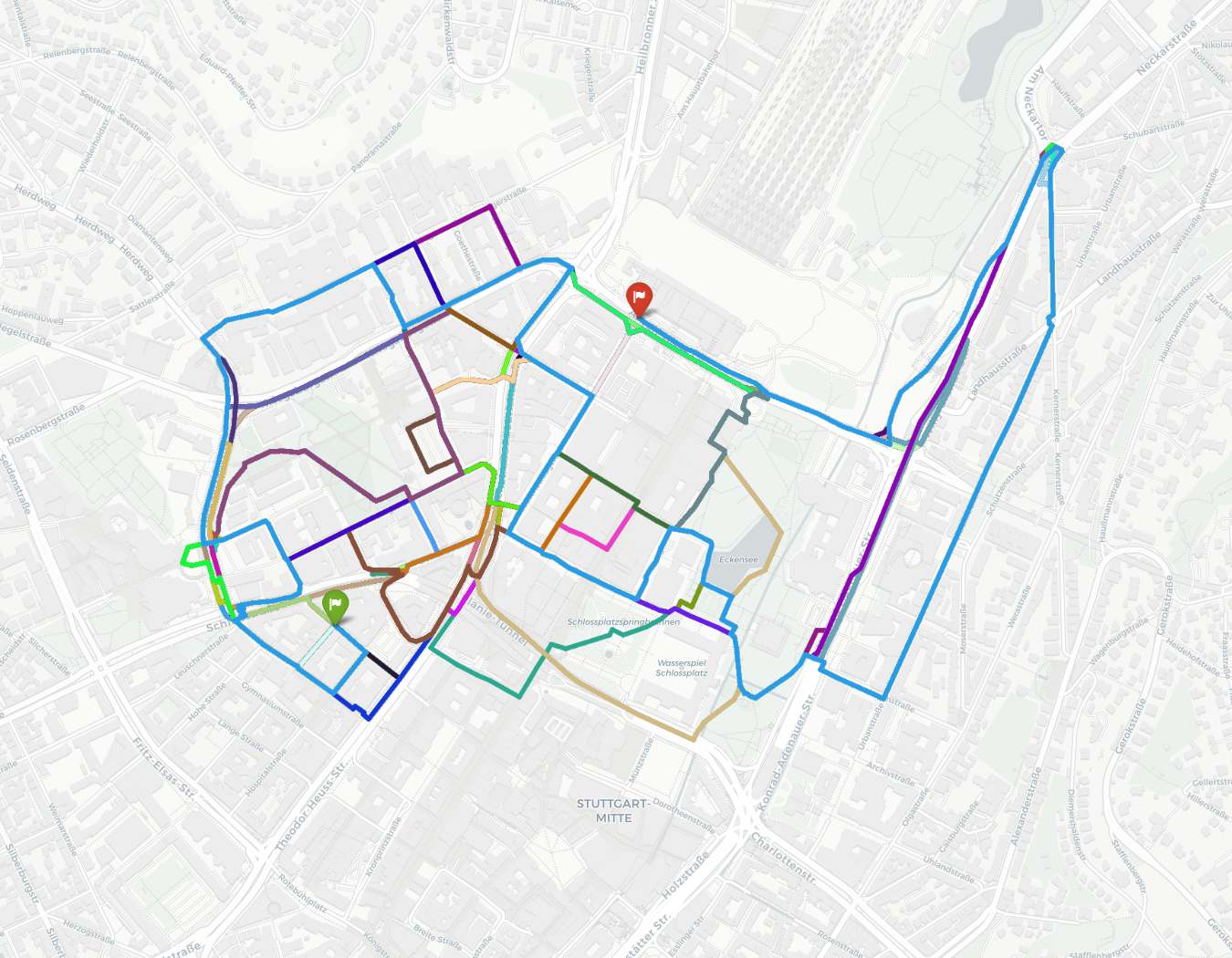}
  
  \caption{800 Meter, 2.155 Routen}
  \label{fig:sub1}
\end{subfigure}%
\begin{subfigure}{.5\textwidth}
  \centering
  \includegraphics[width=0.78\textwidth]{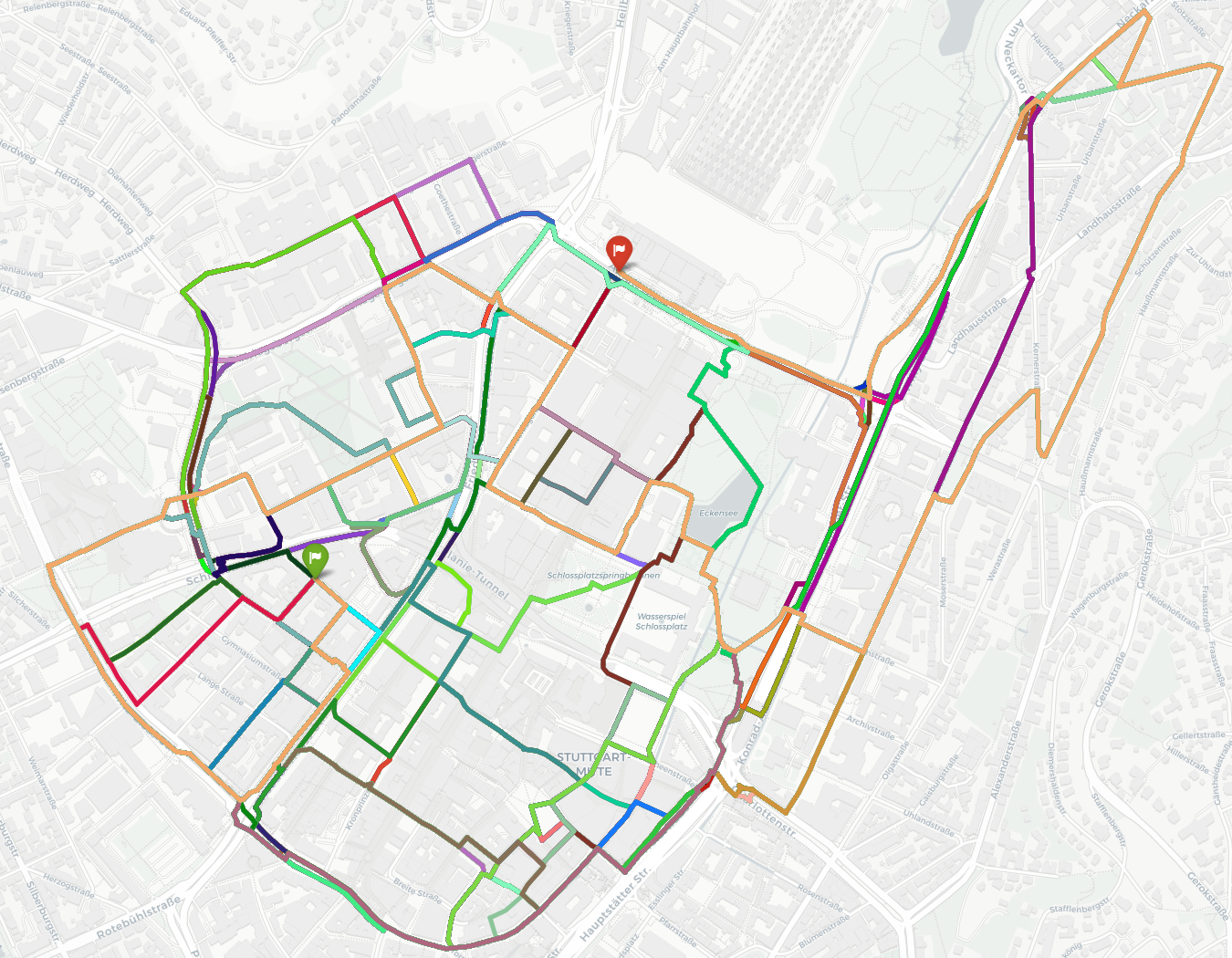}
  
  \caption{1.500 Meter, 2.461 Routen}
  \label{fig:sub2}
\end{subfigure}

\caption{Sichere Routen in Stuttgart auf Graphen verschiedener Größe. (Map data © OpenStreetMap contributors \cite{OpenStreetMap})}
\label{fig:alleWege}
\end{figure}

In Abbildung~\ref{fig:alleWege} sind die resultierenden ordinal effizienten Routen abgebildet. Wird um den Startpunkt ein Quadrat mit 800 Meter halber Seitenlänge gelegt, so besteht der resultierende Graph aus 1099 Knoten und 2824 Kanten. Es werden 2.155 Routen gefunden mit einer durchschnittlichen Länge von 3.219,52 Meter. Hat das Quadrat um den Startpunkt eine halbe Seitenlänge von 1.500 Metern so besteht der resultierende Graph aus 3.097 Knoten und 8.118 Kanten. Es werden 2.461 Routen gefunden mit einer durchschnittlichen Länge von 2.921,99 Metern. Die Anzahl der Routen hat sich erhöht, weil längere Wege mit mehr Kanten in besseren Kategorien gefunden wurden.

\FloatBarrier

\begin{figure}[!htbp]
\centering
\begin{subfigure}{.5\textwidth}
  \centering
  \includegraphics[width=0.78\textwidth]{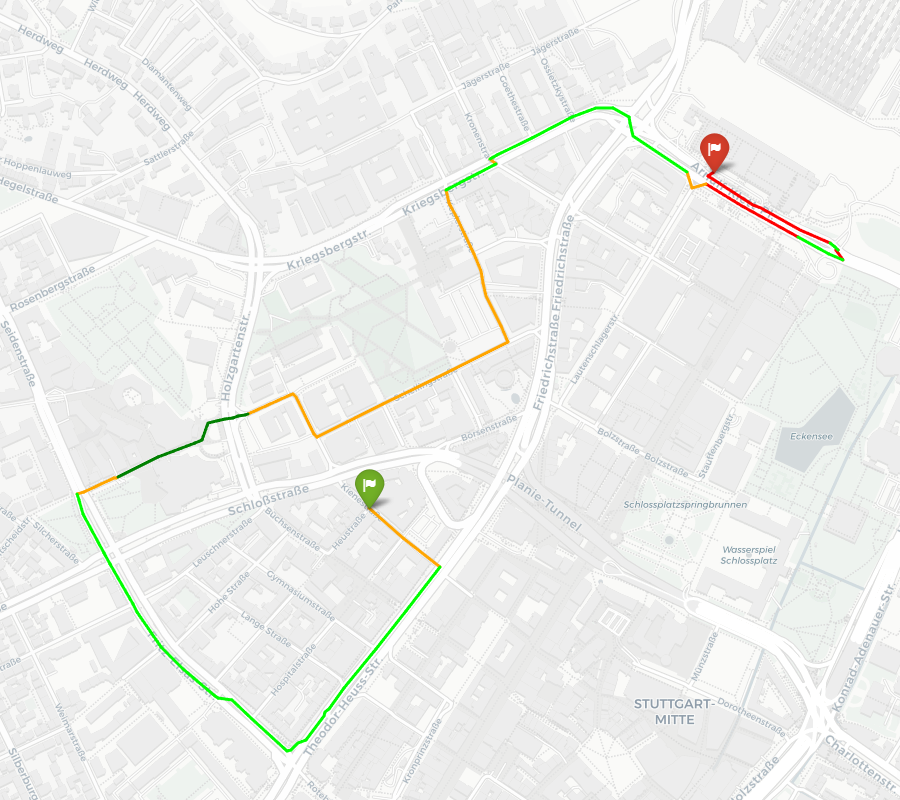}
  
  \caption{Route $\hat{x}$}
  \label{fig:sub1}
\end{subfigure}%
\begin{subfigure}{.5\textwidth}
  \centering
  \includegraphics[width=0.78\textwidth]{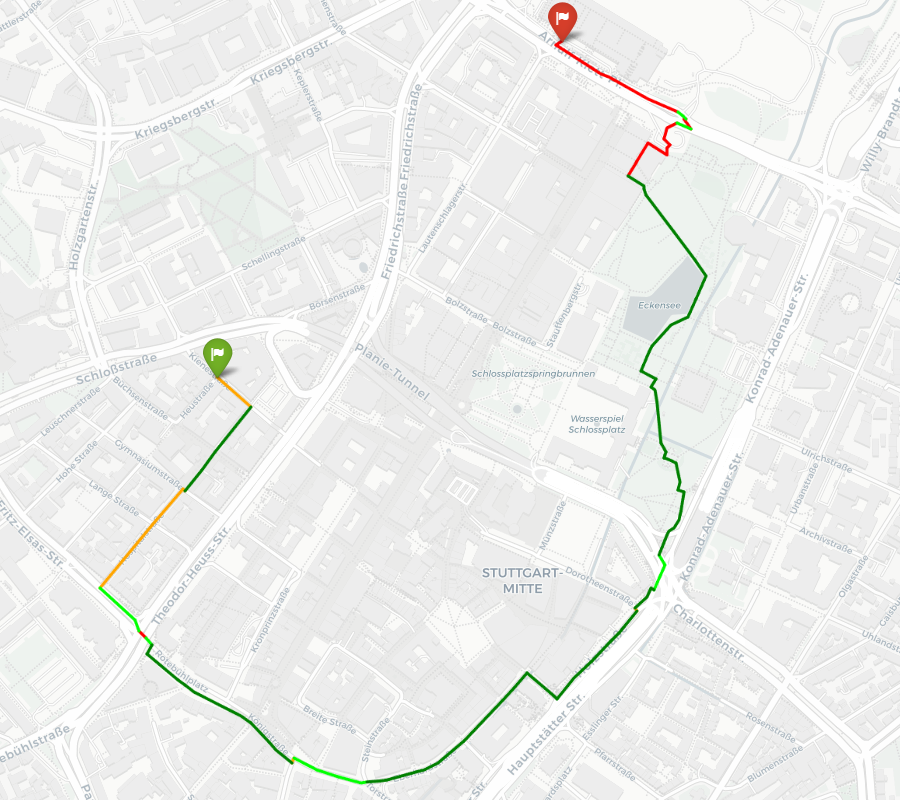}
  
  \caption{Route $\bar{x}$}
  \label{fig:sub2}
\end{subfigure}

\caption{Die Route $\hat{x}$ ist 3.163,63 Meter lang und Route $\bar{x}$ ist 2.873,47 Meter lang. (Map data © OpenStreetMap contributors \cite{OpenStreetMap})}
\label{fig:BeispielRouten}
\end{figure}

Der Unterschied in der durchschnittlichen Routenlänge ist überraschend, weil man intuitiv erwartet, dass ein größerer Graph zu einer größeren durchschnittlichen Routenlänge führt. Dies ist hier jedoch nicht der Fall, da der größere Graph Routen zulässt, die außerhalb des kleineren Graphen verlaufen und dennoch kürzer als die durchschnittliche Routenlänge von 3.219,52 Metern sind, und genau solche Routen sind vermehrt hinzugekommen. Abbildung~\ref{fig:BeispielRouten} zeigt zwei Beispiele $\hat{x}$ und $\bar{x}$ für solche Routen.

Die hohe Anzahl von Routen lässt sich vermutlich darauf zurückführen, dass beide Zielorte sehr zentral gelegen sind und das Straßennetz in Stuttgart äußerst dicht ist. Ähnliche Berechnungen wurden auch für Wuppertal durchgeführt, wo die Anzahl der Routen signifikant niedriger ausfiel (weniger als ein drittel der oben berechneten Routen). Interessanterweise war dabei die Luftlinien-Entfernung zwischen den Start- und Endpunkten ungefähr 300 Meter größer  als in diesem Beispiel.

\begin{figure}[!htbp]
\centering
\begin{subfigure}{.5\textwidth}
  \centering
  \includegraphics[width=0.87\textwidth]{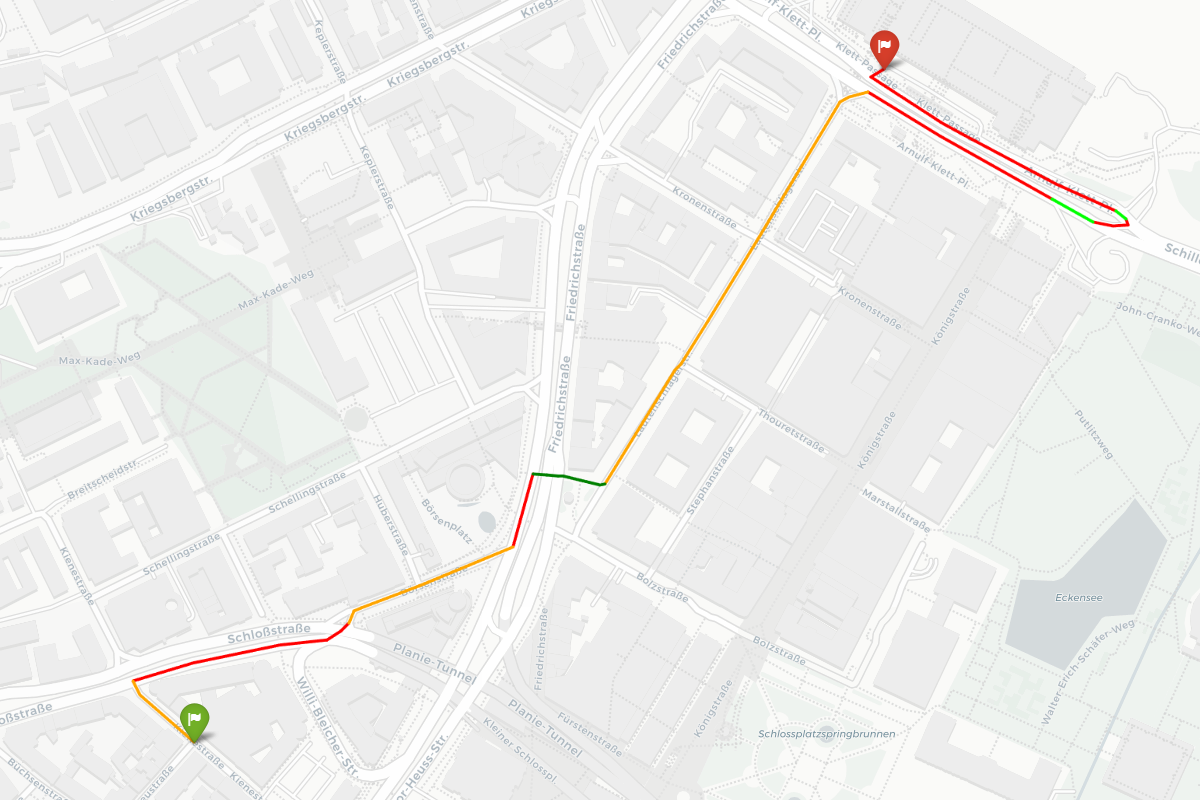}
  
  \caption{Kürzeste Route $x^k$}
  \label{fig:kurerWeg-ungew}
\end{subfigure}%
\begin{subfigure}{.5\textwidth}
  \centering
  \includegraphics[width=0.87\textwidth]{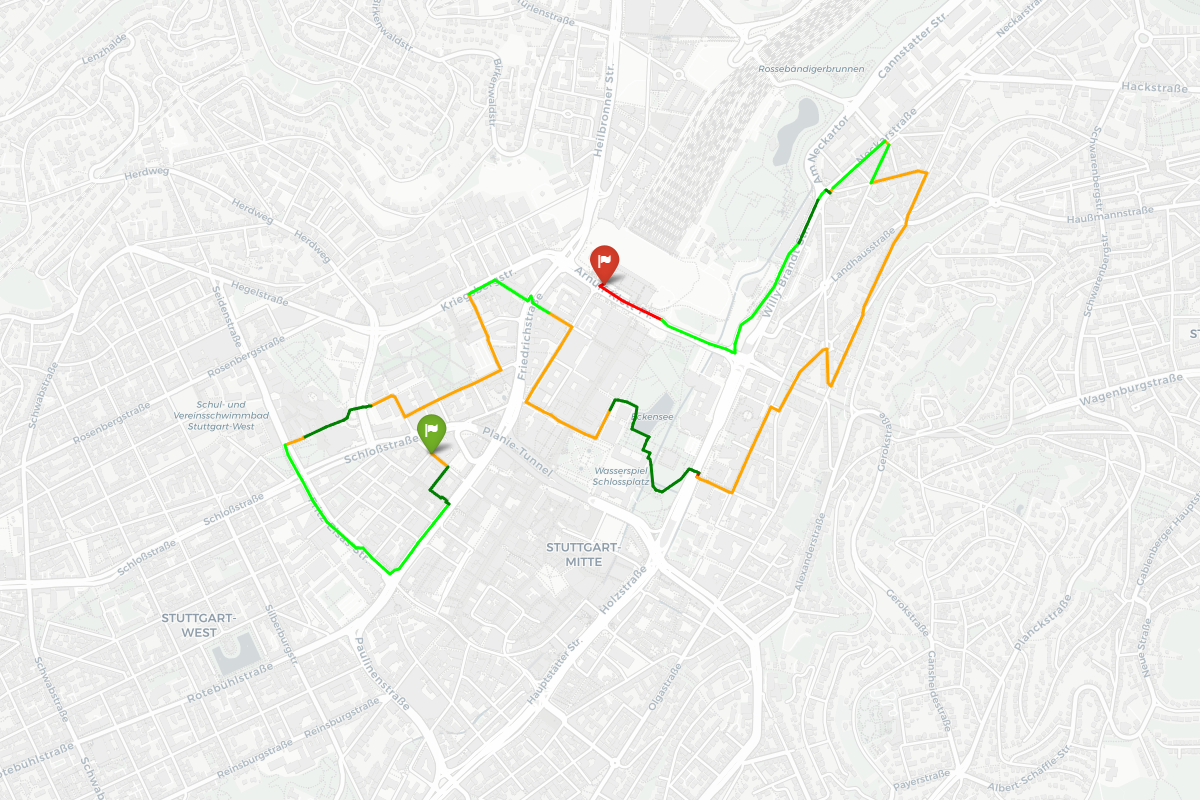}
  
  \caption{Längste Route $x^l$}
  \label{fig:sub2}
\end{subfigure}

\caption{Bilder der Routen $x^k$ und $x^l$ aus dem Graphen der Größe 1.500 Meter. Die einzelnen Abschnitte der Routen sind in der Farbe der jeweiligen Kategorie eingefärbt. (Map data © OpenStreetMap contributors \cite{OpenStreetMap})}
\label{fig:zweiWege}
\end{figure}

Im Folgenden werden die kürzeste Route $x^k$ und die längste Route $x^l$ aus dem Graphen mit einer halben Seitenlänge von 1.500 Metern genauer untersucht. Die Routen sind in Abbildung~\ref{fig:zweiWege} dargestellt.

\FloatBarrier

Die Lösungsvektoren der Routen sind folgende:
\begin{align*}
d(x^k) = \left( \begin{array}{c} 1.364{,}87 \\ 1.307{,}02\\ 1.260{,}55\\ 671{,}40 \end{array} \right),\  \; 
d(x^l) = \left( \begin{array}{c} 7.003{,}06 \\ 5.828{,}10\\ 3.623{,}52\\ 247{,}31 \end{array} \right) \;
\end{align*}

Der erste Eintrag der Kostenvektoren beschreibt die insgesamt zurückgelegte Strecke, der zweite Eintrag beschreibt die zurückgelegte Strecke über hellgrüne, orange und rote Kanten, der dritte Eintrag beschreibt die zurückgelegte Strecke über orange und rote Kanten und der letzte Eintrag beschreibt die zurückgelegte Strecke über rote Kanten.

Die kürzeste Route $x^k$ legt nur 57,85 Meter auf Kanten der besten Kategorie und 46,47 Meter auf Kanten der zweitbesten Kategorie zurück, die restlichen 1.260,55 Meter werden also auf Kanten der dritten und vierten Kategorie zurückgelegt, die deutlich unsicherer sind als die ersten beiden Kategorien. Somit kann die Route $x^k$ als ziemlich unsicher aber sehr kurz deklariert werden.

Route $x^l$ ist die längste berechnete Route und ist etwa fünfmal so lang wie Route $x^k$. Dafür fährt man mit der Route $x^l$ hauptsächlich über Kanten, die nicht in der schlechtesten Kategorie sind, denn Route $x^l$ nutzt nur 247,31 Meter rote Kanten. Allerdings ist die Gesamtstrecke über orange und rote Kanten schon etwa dreimal solang bei Route $x^l$ verglichen mit Route $x^k$.

Offensichtlich ist die kürzeste Route $x^k$ für Radfahrer denen Sicherheit wichtig ist, nicht die beste Wahl, insbesondere wenn man zusammen mit Kindern unterwegs ist, weil fast der gesamte Weg über Straßen ohne Fahrradwege verläuft (orange und rote Kanten). Daraus folgt, dass es sinnvoll ist die Gewichte so anzupassen, dass man kleine Umwege in Kauf nimmt um unsichere Straßen zu vermeiden. 

Da insbesondere der Unterschied zwischen Kanten mit und ohne Fahrradweg, also zwischen hellgrünen und orangen Kanten, besonders groß ist, werden im Folgenden die Gewichte $\omega_1= \omega_3 = 1{,}5$ und $\omega_2 = 2{,}5$ verwendet. Dies bedeutet, dass die 1,5-fache Entfernung auf grünen Straßen gegenüber hellgrünen Straßen bevorzugt wird, ebenso die 1,5-fache Entfernung auf orangen Straßen gegenüber roten Straßen und die 2,5-fache Entfernung auf hellgrünen Straßen gegenüber orangen Straßen.

\begin{figure}[!htbp]
\centering
\begin{subfigure}{.5\textwidth}
  \centering
  \includegraphics[width=0.8\textwidth]{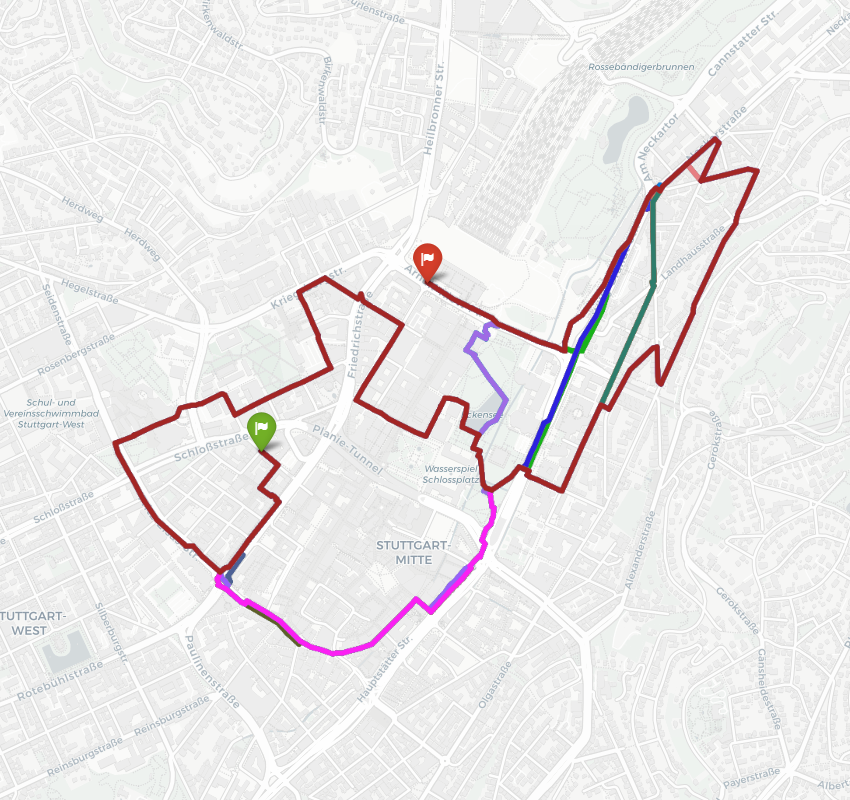}
  
  \caption{1.500 Meter, 57 Routen}
  \label{fig:alleGewWege}
\end{subfigure}%
\begin{subfigure}{.5\textwidth}
  \centering
\includegraphics[width=0.8\textwidth]{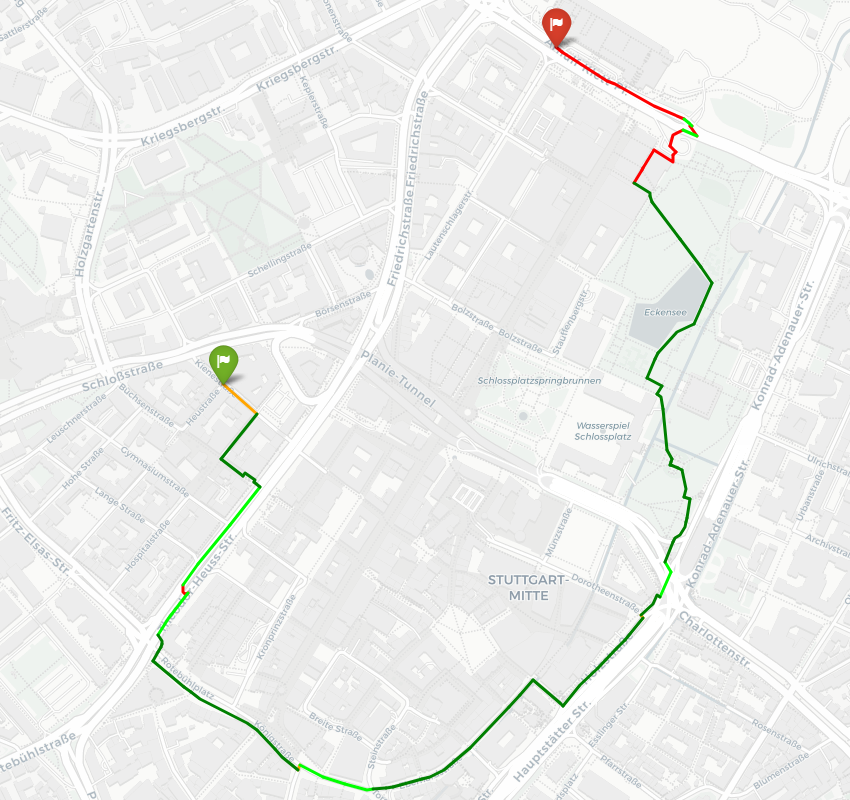}
  
  \caption{Kürzeste gewichtete Route $x^g$}
  \label{fig:kurz-gewichteteWege}
\end{subfigure}

\caption{Bilder der Routen die mit $\omega_1 = \omega_3 = 1{,}5$ und $\omega_2 = 2{,}5$ berechnet wurden. (Map data © OpenStreetMap contributors \cite{OpenStreetMap})}
\label{fig:AlleGewichteRouten}
\end{figure}

Mit diesen Gewichten werden auf dem Graph der Größe 1.500 nur noch 57 optimale Wege mit einer durchschnittlichen Länge von 4.465,20 Metern ermittelt, welche in Abbildung~\ref{fig:alleGewWege} dargestellt sind.  Da die längste Route unverändert bleibt, werden wir uns im folgenden nur die kürzeste Route genauer anschauen. Diese ist in Abbildung~\ref{fig:kurz-gewichteteWege} dargestellt.

Der gewichteten Lösungsvektoren von  $x^g$ sieht wie folgt aus:
\begin{align*}
d^{(\omega)}(x^g) = \left( \begin{array}{c} 5.387{,}30 \\ 2.353{,}95\\ 741{,}88\\ 437{,}26 \end{array} \right),\  \; 
\end{align*}

Da der gewichteten Lösungsvektor nicht so gut untersucht werden kann wie der ungewichteten, ist hier der ungewichteten Lösungsvektor der Routen $x^g$ dargestellt:
\begin{align*}
d(x^g) = \left( \begin{array}{c} 2.878{,}88 \\ 1.022{,}51 \\ 523{,}25 \\ 437{,}26 \end{array} \right),\  \; 
\end{align*}

Die Distanz der kürzesten Route $x^g$ der gewichteten Routen hat sich im Vergleich zur kürzesten ungewichteten Route $x'$ etwa verdoppelt. Abbildung~\ref{fig:kurz-gewichteteWege} und Abbildung~\ref{fig:kurerWeg-ungew} zeigen, dass sich die Routen deutlich voneinander unterscheiden. Insbesondere fällt auf, dass die Route $x^g$ überwiegend auf Radwegen verläuft, was bei der Route $x^k$ nicht der Fall ist. Vergleicht man nun die Routen anhand der Vektoren, so zeigt sich, dass die Route $x^g$ um 1.514,01 Meter länger ist als die kürzeste Route $x^k$, dafür aber 737,30 Meter weniger auf orangen und roten Straßen zurückgelegt werden. Außerdem werden auf der Route $x^g$ 2.355,63 Meter auf Radwegen zurückgelegt, was 81,8\% der Gesamtstrecke entspricht. Aufgrund dieser Beobachtungen kann man sagen, dass die Route $x^g$ eine bevorzugte Alternative zur kürzesten Route $x^k$ ist.

\begin{table}[!htbp]
\begin{tabular}[h]{l|c|c|c|r}
$\omega_i$ &  Anzahl der Routen & Durchschnittliche Routen Länge [m] & Berechnungszeit [sec] \\[0.5ex]
\hline
$16$ &  $4$ & $5.635{,}72$ & $0{,}35$ \\[0.5ex]
$8$ &  $8$ & $4.933{,}45$ & $0{,}53$ \\[0.5ex]
$4$ & $19$ & $4.806{,}27$ & $0{,}88$ \\[0.5ex]
$2$ &  $37$ & $4.313{,}67$ & $2{,}32$ \\[0.5ex]
$1$ &  $2.461$ & $2.921{,}99$ & $2.207{,}58$ \\[0.5ex]
\end{tabular}

\caption{Diese Tabelle zeigt die Anzahl der Routen und ihre durchschnittliche Länge in
Meter für verschiedene Gewichtungen $\omega$. Die Berechnungszeit ist in Sekunden angegeben. Die Routen wurden auf einem
Graphen mit einer Größe von 1.500 Metern berechnet. Die Werte für $\omega_i$  in der Tabelle stellen alle Einträge in den Vektoren $\omega$ dar.}
\label{tab:gewichte}
\end{table}

Abschließend haben wir den Einfluss der Gewichte auf die Anzahl der Lösungen und die durchschnittliche Routenlänge untersucht. Die Ergebnisse sind in Tabelle~\ref{tab:gewichte} dargestellt. Je größer $\omega$ wird,
desto weniger Routen werden gefunden und desto geringer ist die Rechenzeit. Die Tabelle zeigt, dass bei einer Erhöhung von $\omega$ (von 1 aus) die kürzeren Routen herausgefiltert werden. 

Eine gute Wahl der Gewichte ist daher notwendig, um das volle Potenzial der ordinalen kürzesten Wege auszuschöpfen und eine geeignete Anzahl von Routen zu erhalten, die für Radfahrer ansprechend sind. Man sollte bedenken, dass es nicht immer optimal ist $\omega$ 
möglichst groß zu wählen, auch wenn dies die Rechenzeit und die Anzahl der Routen reduziert. Denn eine Auswahl an Routen ist in der Praxis sinnvoll, da die präferierte Route je nach Radfahrer und Situation variieren kann. Erfahrene Radfahrer trauen es sich vielleicht eher zu auch gefährliche Routen zu wählen und der Zeitdruck beeinflusst wie wichtig einem die Länge der Route ist. Denn der große Vorteil der Berechnung ordinal kürzester Wege ist, dass man mehrere Routen erhält und diese nach möglichen Bedürfnissen und Kriterien filtern kann, um somit eine Route zu finden, die einem Radfahrer zusagt. Nichtsdestotrotz ist es sinnvoll $\omega$ größer als 1 zu wählen, weil die Kategorien so klarer von einander getrennt werden und Wege auf sichereren Straßen stärker bevorzugt werden. 

\section{Zusammenfassung und Ausblick}
In diesem Beitrag wurde eine Möglichkeit zur Modellierung ungewichteter und gewichteter ordinaler Kosten vorgestellt. Diese Modellierung hat den Vorteil, dass bekannte Verfahren zur Lösung multikriterieller Optimierungsprobleme genutzt werden können, um Probleme mit (gewichteten) ordinalen Kosten zu lösen. Dies wurde genutzt um sichere Fahrradwege in Stuttgart auf Basis von OpenStreetMaps Daten und mit Hilfe des Multikriteriellen Dijkstra Algorithmus zu ermitteln. Dabei wurde der Nutzen von gewichteten ordinalen Kosten deutlich, denn die Gewichte sorgen für kürze Rechenzeiten und eine Beschränkung der Lösungen auf praxistaugliche Routen, da sehr kurze aber unsichere Routen eliminiert werden. Die Gewichte sorgen für eine stärkere Abgrenzung zwischen den Kategorien und sind immer dann besonders nützlich, wenn man bereit ist Umwege für die Sicherheit in Kauf zu nehmen, also beispielsweise wenn man mit Kindern unterwegs ist oder eine Radtour ohne Zeitdruck unternehmen möchte. Dennoch erhält man meist mehrere Routen als mögliche Lösungen, was dem einzelnen Radfahrer erlaubt nach persönlicher Präferenz die für ihn passende zu wählen. 

Um die Auswahl der Routen noch besser an die individuellen Wünsche der Radfahrer anpassen zu können, wäre es interessant zusätzliche Gewichte zu berücksichtigen, um zu lange Routen zu vermeiden. Diese könnten analog zu den $\omega$-Gewichten beschreiben, wie lang ein Umweg maximal sein darf um ein Stück unsicheren Weg zu vermeiden. Es wäre spannend zu untersuchen, wie sich die Lösungen in unserem Anwendungsfall dann verändern.

Darüber hinaus könnte man in zukünftiger Forschung, versuchen die Modellierung weiter zu verbessern.
Die aktuelle Modellierung kategorisiert Sicherheit primär auf Grund der Infrastruktur der Straßen. In der Praxis wird die Sicherheit beispielsweise auch von der Dichte und Geschwindigkeit des Verkehrs, der Beleuchtung bei Dunkelheit, der Steigung und vielen weiteren Faktoren beeinflusst. Um das Modell zu verbessern sollten diese Faktoren bei der Kategorisierung der Straßen mit berücksichtigt werden. Des weiteren sind insbesondere Kreuzungen eine potentielle Gefahr für Radfahrer. Diese sind bisher unzureichend modelliert und haben, auf Grund der geringen Strecke die beim Überqueren einer Kreuzung zurück gelegt werden muss, im aktuellen Modell einen geringen Einfluss. Zukünftige Forschung sollte sich mit der Frage befassen, wie Kreuzungen besser modelliert werden können, sodass auch die Überquerungsrichtung mit einbezogen wird, denn im Allgemeinen ist das Linksabbiegen für einen Radfahrer wesentlich gefährlicher als das Rechtsabbiegen. Dadurch könnte die Qualität der berechneten Routen weiter verbessert werden.

\bibliography{literatur}

\begin{thebibliography}{27}
\providecommand{\natexlab}[1]{#1}
\providecommand{\url}[1]{\texttt{#1}}
\expandafter\ifx\csname urlstyle\endcsname\relax
  \providecommand{\doi}[1]{doi: #1}\else
  \providecommand{\doi}{doi: \begingroup \urlstyle{rm}\Url}\fi

\bibitem[{British Olympic Council}(1909)]{Olympia1908}
{British Olympic Council}.
\newblock Fourth olympiad; being the official report of the olympic games of
  1908 celebrated in london under the patronage of his most gracious majesty
  king edward vii and by the sanction of the international olympic committee.
\newblock Digitally published by the LA84 Foundation, March 1909.
\newblock URL
  \url{https://digital.la84.org/digital/collection/p17103coll8/id/8217/rec/6}.

\bibitem[{International Olympic Committee}(1913)]{Olympia2012}
{International Olympic Committee}.
\newblock Fifth olympiad: the official report of the olympic games of
  stockholm, 1912 swedish olympic committee.
\newblock Digitally published by the LA84 Foundation, 1913.
\newblock URL
  \url{https://digital.la84.org/digital/collection/p17103coll8/id/11660/rec/7#page=1178}.

\bibitem[Sitarz(2013)]{Sitarz2013medal}
Sebastian Sitarz.
\newblock The medal points' incenter for rankings in sport.
\newblock \emph{Applied Mathematics Letters}, 26\penalty0 (4):\penalty0
  408--412, apr 2013.
\newblock \doi{10.1016/j.aml.2012.10.014}.

\bibitem[Gomes~Júnior et~al.(2014)Gomes~Júnior, Mello, and
  Angulo-Meza]{GomesJunior2014Sequential}
Silvio Gomes~Júnior, João Mello, and Lidia Angulo-Meza.
\newblock Sequential use of ordinal multicriteria methods to obtain a ranking
  for the 2012 summer olympic games.
\newblock \emph{WSEAS Transactions on Systems}, 13:\penalty0 223--230, 04 2014.

\bibitem[Du(2018)]{Du2018Modifying}
Jiangze Du.
\newblock Modifying olympics medal table via a stochastic multicriteria
  acceptability analysis.
\newblock \emph{Mathematical Problems in Engineering}, 2018:\penalty0 1--11,
  aug 2018.
\newblock \doi{10.1155/2018/8729158}.

\bibitem[Perini et~al.(2022)Perini, Langville, Kramer, Shrager, and
  Shapiro]{Perini2022Weight}
Tyler Perini, Amy Langville, Glenn Kramer, Jeff Shrager, and Mark Shapiro.
\newblock Weight set decomposition for weighted rank aggregation: An
  interpretable and visual decision support tool, 2022.
\newblock URL \url{https://arxiv.org/abs/2206.00001}.

\bibitem[Bartee(1971)]{Bartee1971Problem}
Edwin~M. Bartee.
\newblock Problem solving with ordinal measurement.
\newblock \emph{Management Science}, 17\penalty0 (10):\penalty0 B622--B633,
  1971.
\newblock ISSN 00251909, 15265501.
\newblock URL \url{http://www.jstor.org/stable/2628998}.

\bibitem[Fishburn and LaValle(1996)]{Fishburn1996Binary}
Peter~C. Fishburn and Irving~H. LaValle.
\newblock Binary interactions and subset choice.
\newblock \emph{European Journal of Operational Research}, 92\penalty0
  (1):\penalty0 182--192, 1996.
\newblock ISSN 0377-2217.
\newblock \doi{10.1016/0377-2217(95)00073-9}.

\bibitem[Brams et~al.(2003)Brams, Edelman, and Fishburn]{Brams2003Fair}
Steven~J Brams, Paul~H Edelman, and Peter~C Fishburn.
\newblock Fair division of indivisible items.
\newblock \emph{Theory and Decision}, 55\penalty0 (2):\penalty0 147--180, 2003.

\bibitem[Brams and King(2005)]{Brams2005Efficient}
Steven~J. Brams and Daniel~L. King.
\newblock Efficient fair division: Help the worst off or avoid envy?
\newblock \emph{Rationality and Society}, 17\penalty0 (4):\penalty0 387--421,
  2005.
\newblock \doi{10.1177/1043463105058317}.

\bibitem[Bouveret and Endriss(2010)]{Bouveret2010Fair}
Sylvain Bouveret and Ulle Endriss.
\newblock Fair division under ordinal preferences: Computing envy-free
  allocations of indivisible goods.
\newblock In \emph{Proceedings of the 19th European Conference on Artificial
  Intelligence (ECAI-2010)}, volume 215, pages 387--392, 06 2010.
\newblock \doi{10.3233/978-1-60750-606-5-387}.

\bibitem[Schäfer et~al.(2020)Schäfer, Dietz, Fröhlich, Ruzika, and
  Figueira]{Schafer:SP}
Luca~E. Schäfer, Tobias Dietz, Nicolas Fröhlich, Stefan Ruzika, and José~R.
  Figueira.
\newblock Shortest paths with ordinal weights.
\newblock \emph{European Journal of Operational Research}, 280\penalty0
  (3):\penalty0 1160–1170, 2020.
\newblock ISSN 0377-2217.
\newblock \doi{10.1016/j.ejor.2019.08.008}.

\bibitem[Klamroth et~al.(2023{\natexlab{a}})Klamroth, Stiglmayr, and
  Sudhoff]{Klamroth2023Multi}
K.~Klamroth, M.~Stiglmayr, and J.~Sudhoff.
\newblock Multi-objective matroid optimization with ordinal weights.
\newblock \emph{Discrete Applied Mathematics}, 335:\penalty0 104--119,
  2023{\natexlab{a}}.
\newblock \doi{10.1016/j.dam.2022.07.017}.

\bibitem[O'Mahony and Wilson(2013)]{OMahony2013Sorted}
Conor O'Mahony and Nic Wilson.
\newblock Sorted-pareto dominance and qualitative notions of optimality.
\newblock In Linda~C. van~der Gaag, editor, \emph{Symbolic and Quantitative
  Approaches to Reasoning with Uncertainty}, pages 449--460. Springer Berlin,
  Heidelberg, 2013.
\newblock ISBN 978-3-642-39091-3.

\bibitem[Schweigert(1999)]{Schweigert1999Ordered}
Dietmar Schweigert.
\newblock Ordered graphs and minimal spanning trees.
\newblock \emph{Foundations of computing and decision sciences}, 24\penalty0
  (2):\penalty0 219--229, 1999.

\bibitem[Bossong and Schweigert(1999)]{Bossong1999Minimal}
Ulrike Bossong and Dietmar Schweigert.
\newblock Minimal paths on ordered graphs.
\newblock 1999.
\newblock URL \url{http://nbn-resolving.de/urn:nbn:de:hbz:386-kluedo-4666}.

\bibitem[Delort et~al.(2011)Delort, Spanjaard, and Weng]{Delort2011Committee}
Charles Delort, Olivier Spanjaard, and Paul Weng.
\newblock Committee selection with a weight constraint based on a pairwise
  dominance relation.
\newblock In Ronen~I. Brafman, Fred~S. Roberts, and Alexis Tsouki{\`a}s,
  editors, \emph{Algorithmic Decision Theory}, pages 28--41, Berlin,
  Heidelberg, 2011. Springer Berlin Heidelberg.
\newblock ISBN 978-3-642-24873-3.
\newblock \doi{10.1007/978-3-642-24873-3_3}.

\bibitem[Schäfer et~al.(2021)Schäfer, Dietz, Barbati, Figueira, Greco, and
  Ruzika]{SCHAFER:knapsack}
Luca~E. Schäfer, Tobias Dietz, Maria Barbati, José~Rui Figueira, Salvatore
  Greco, and Stefan Ruzika.
\newblock The binary knapsack problem with qualitative levels.
\newblock \emph{European Journal of Operational Research}, 289\penalty0
  (2):\penalty0 508–514, 2021.
\newblock \doi{10.1016/j.ejor.2020.07.040}.

\bibitem[Klamroth et~al.(2023{\natexlab{b}})Klamroth, Stiglmayr, and
  Sudhoff]{Klamroth2023Ordinal}
K.~Klamroth, M.~Stiglmayr, and J.~Sudhoff.
\newblock Ordinal optimization through multi-objective reformulation.
\newblock \emph{European Journal of Operational Research}, 311\penalty0
  (2):\penalty0 427--443, 2023{\natexlab{b}}.
\newblock \doi{10.1016/j.ejor.2023.04.042}.

\bibitem[Klamroth et~al.(2023{\natexlab{c}})Klamroth, Stiglmayr, and
  Santos]{Klamroth2023Weighted}
K.~Klamroth, M.~Stiglmayr, and J.~Sudhoff Santos.
\newblock Weighted ordinal cones: Modelling preferences in ordinal
  optimization.
\newblock wird bis Januar auf arXiv verfügbar sein, 2023{\natexlab{c}}.

\bibitem[geovelo()]{geovelo}
geovelo.
\newblock Geovelo, an application for short and safe bicycle route
  computations.
\newblock \url{https://geovelo.fr}, 2022.
\newblock Accessed: 2022-01-28.

\bibitem[Kergosien et~al.(2021)Kergosien, Giret, N{\'e}ron, and
  Sauvanet]{kerg:anef:2021}
Yannick Kergosien, Antoine Giret, Emmanuel N{\'e}ron, and Ga{\"e}l Sauvanet.
\newblock An efficient label-correcting algorithm for the multiobjective
  shortest path problem.
\newblock \emph{INFORMS Journal on Computing}, 34\penalty0 (1):\penalty0
  76--92, 2021.
\newblock \doi{10.1287/ijoc.2021.1081}.

\bibitem[Sauvanet and N{\'e}ron(2010)]{sauv:sear:2010}
G.~Sauvanet and E.~N{\'e}ron.
\newblock Search for the best compromise solution on multiobjective shortest
  path problem.
\newblock \emph{Electron. Notes Discrete Math.}, 36:\penalty0 615--622, 2010.
\newblock \doi{10.1016/j.endm.2010.05.078}.

\bibitem[Ehrgott(2005)]{Ehrgott2005Multicriteria}
M.~Ehrgott.
\newblock \emph{Multicriteria Optimization}.
\newblock Springer Verlag, Berlin,~Heidelberg, 2005.
\newblock \doi{10.1007/3-540-27659-9}.

\bibitem[de~las Casas et~al.(2021)de~las Casas, Sede{\~{n}}o-Noda, and
  Borndörfer]{Casas2021Improved}
P.~M. de~las Casas, A.~Sede{\~{n}}o-Noda, and R.~Borndörfer.
\newblock An improved multiobjective shortest path algorithm.
\newblock \emph{Computers {\&} Operations Research}, 135:\penalty0 105424, nov
  2021.
\newblock \doi{10.1016/j.cor.2021.105424}.

\bibitem[{OpenStreetMap contributors}(2017)]{OpenStreetMap}
{OpenStreetMap contributors}.
\newblock {Planet dump retrieved from https://planet.osm.org at 29.09.2023}.
\newblock \url{ https://www.openstreetmap.org }, 2017.
\newblock {Open Data Commons Open Database License (ODbL)
  \url{www.opendatacommons.org/licenses/odbl}}.

\bibitem[Boeing(2017)]{osmnx}
G.~Boeing.
\newblock Osmnx: New methods for acquiring, constructing, analyzing, and
  visualizing complex street networks.
\newblock \emph{Computers, Environment and Urban Systems}, pages 126--139,
  2017.
\newblock \doi{10.1016/j.compenvurbsys.2017.05.004}.

\end{thebibliography}
\end{document}